\documentclass[12pt,a4paper]{amsart}

\usepackage{amssymb}
\usepackage{amsrefs} 
\usepackage[all]{xy}

\allowdisplaybreaks

\newcommand{\bbA}{\mathbb{A}}
\newcommand{\bbD}{\mathbb{D}}
\newcommand{\bbN}{\mathbb{N}}
\newcommand{\bbZ}{\mathbb{Z}}

\newcommand{\fx}{\mathfrak{x}}
\newcommand{\fz}{\mathfrak{z}}
\newcommand{\fP}{\mathfrak{P}}
\newcommand{\fR}{\mathfrak{R}}
\newcommand{\fS}{\mathfrak{S}}
\newcommand{\fT}{\mathfrak{T}}

\let\Im=\undefined
\DeclareMathOperator{\Du}{D} %
\DeclareMathOperator{\Id}{Id} %
\DeclareMathOperator{\Im}{Im} %
\DeclareMathOperator{\Dom}{Dom} %
\DeclareMathOperator{\End}{End} %
\DeclareMathOperator{\Ext}{Ext} %
\DeclareMathOperator{\Hom}{Hom} %
\DeclareMathOperator{\rad}{rad} %
\DeclareMathOperator{\Coker}{Coker} %

\newtheorem{lemm}{Lemma}[section]
\newtheorem{theo}[lemm]{Theorem}

\newcommand{\ol}{\overline}

\numberwithin{equation}{subsection}

\title[A characterization of admissible algebras]%
{A characterization of admissible algebras with formal two-ray
modules}

\author{Grzegorz Bobi\'nski}

\thanks{Supported by the Polish Research Grant KBN 5 PO3A
008 21}

\address{Faculty of Mathematics and Computer Science \\
Nicolaus Coperincus University \\ Chopina 12/18 \\ 87-100 Toru\'n
\\ Poland}

\email{gregbob@mat.uni.torun.pl}



\begin{document}

\begin{abstract}
In the paper we characterize, in terms of quivers and relations,
the admissible algebras with formal two-ray modules introduced
by G.~Bobi\'nski and A.~Skowro\'nski 
[\textit{Cent. Eur. J.Math.}~\textbf{1} (2003), 457--476]. 
\end{abstract}

\maketitle

Throughout the paper $K$ denotes a fixed algebraically closed
field. By an algebra we mean a finite dimensional $K$-algebra with
identity and by a module a finite dimensional (left) module.

According to the Drozd's Tame and Wild Theorem~\cite{Dr} (see
also~\cite{CB2}) the representation-infinite algebras can be
divided into two disjoint classes. One class consists of the wild
algebras whose representation theory comprises the representation
theories of all algebras. The second class consists of the tame
algebras, for which in each dimension all but finitely many
indecomposable modules can be parameterized by a finite number of
lines (see also~\cite{CB3}). Thus one may realistically hope to
classify the indecomposable modules only for the tame algebras.
The first level in the hierarchy of the tame algebras is occupied
by the domestic algebras. These are characterized by the property
that there exists a common bound for the numbers of lines
necessary to classify the indecomposable modules of a given
dimension (see~\cite{Ri1}). The representation theory of all
strongly simply connected domestic (more generally, of polynomial
growth) algebras seems to be well-understood (see for
example~\cites{MaSkTo, PeSk, Sk1, Sk2, SkZw}). For example, if $A$
is a strongly simply connected domestic algebra, then all but
finitely many components of the Auslander--Reiten quiver of $A$
are homogeneous tubes. An important and interesting open problem
is to establish the representation theory of arbitrary domestic
algebras.

In~\cite{BoSk1} (being a continuation of earlier
works~\cites{BoDrSk, BoSk2}) Skowro\'nski and the author
introduced a class of domestic algebras, called admissible
algebras with formal two-ray modules. The Aus\-lan\-der--Reiten
quivers of these algebras revealed new interesting properties of
the Auslander--Reiten quivers of the domestic algebras. It is also
expected that the representation theory of an arbitrary domestic
algebra is approximated by the representation theory of an
admissible algebra with formal two-ray modules. A drawback of the
construction given in~\cite{BoSk1} is that it does not give a
handy criterion to determine if a given algebra belongs to the
considered class. In this paper we characterize the above algebras
in terms of their quivers and relations. We also translate to this
new language, the description of the Auslander--Reiten quiver of
these algebras.

For basic background on the representation theory of algebras and
all unexplained notions we refer to~~\cite{Ri2}.

The paper is organized as follows. In Section~\ref{sectmain} we
present the main results of the paper, in Section~\ref{sectcomb}
we develop necessary combinatorics, while in final
Section~\ref{sectproof} we prove the main theorems.

\section{Main results} \label{sectmain}

The aim of this section is to present the main results of the
paper. First we introduce a necessary notation and definitions. In
the paper, by $\bbN$ (respectively, $\bbN_0$) we denote the set of
positive (nonnegative) integers. If $m$ and $n$ are integers, then
$[m, n]$ denotes the set of all integers $k$ such that $m \leq k
\leq n$. For a sequence $f : [1, n] \to \bbN$, $n \in \bbN_0$, of
positive integers, we denote $n$ by $|f|$. We identify the finite
subsets of $\bbN$ with the corresponding increasing sequences of
positive integers. In particular, if $F$ is a finite subset of
$\bbN$ and $i \in [1, |F|]$, then $F_i$ denotes the $i$-th element
of $F$ with respect to the usual ordering of integers.

By a defining system we mean a quadruple $(p, q, S, T)$, where $p$
and $q$ are sequences of positive integers with $|q| = |p|$ and
$\sum_{i = 1}^{|p|} p_i \geq 2$, and $S = (S_i)_{i = 1}^{|p|}$ and
$T = (T_i)_{i = 1}^{|p|}$ are families of subsets of $\bbN$ such
that for each $i \in [1, |p|]$ the following conditions hold: $T_i
\subseteq S_i \subseteq [2, p_i + |T_i|]$, $p_i + |T_i| \not \in
T_i$, and $j + 1 \not \in S_i$ for $j \in S_i$. We write $T_{i,
j}$ instead of $(T_i)_j$ for $i \in [1, |p|]$ and $j \in [1,
|T_i|]$.

For a defining system $(p, q, S, T)$ we define a quiver $Q$ in the
following way: the vertices of $Q$ are
\begin{itemize}

\item
$x_{i, j}$, $i \in [1, |p|]$, $j \in [0, p_i + |T_i|]$,

\item
$y_{i, j}$, $i \in [1, |p|]$, $j \in [1, q_i - 1]$,

\item
$z_{i, j}$, $i \in [1, |p|]$, $j \in S_i$,

\end{itemize}
and the arrows of $Q$ are
\begin{itemize}

\item
$\alpha_{i, j} : x_{i, j} \rightarrow x_{i, j - 1}$, $i \in [1,
|p|]$, $j \in [1, p_i + |T_i|]$,

\item
$\beta_{i, j} : y_{i, j} \rightarrow y_{i, j - 1}$, $i \in [1,
|p|]$, $j \in [1, q_i]$, where $y_{i, 0} = x_{i + 1, 0}$ (with
$x_{n + 1, 0} = x_{1, 0}$) and $y_{i, q_i} = x_{i, p_i}$,

\item
$\gamma_{i, j} : z_{i, j} \rightarrow x_{i, j}$, $i \in [1, |p|]$,
$j \in S_i$,

\item
$\xi_{i, j} : x_{i, p_i + j} \rightarrow z_{i, T_{i, j}}$, $i \in
[1, |p|]$, $j \in [1, |T_i|]$.

\end{itemize}
Let $A$ be the path algebra of the quiver $Q$ bounded by
relations:
\begin{itemize}

\item
$\alpha_{i, j - 1} \alpha_{i, j} \gamma_{i, j}$, $i \in [1, |p|]$,
$j \in S_i$,

\item
$\beta_{i, q_i} \alpha_{i, p_i + 1}$, $i \in [1, |p|]$ such that
$|T_i| > 0$,

\item
$\xi_{i, j - 1} \alpha_{i, p_i + j}$, $i \in [1, |p|]$, $j \in [2,
|T_i|]$,

\item
$\alpha_{i, T_{i, j}} \gamma_{i, T_{i, j}} \xi_{i, j} - \alpha_{i,
T_{i, j}} \alpha_{i, T_{i, j} + 1} \cdots \alpha_{i, p_i + j - 1}
\alpha_{i, p_i + j}$, $i \in [1, |p|]$, $j \in [1, |T_i|]$.

\end{itemize}
We say that $Q$ is the quiver and $A$ is the algebra of the
defining system $(p, q, S, T)$. For example, if $p = (6, 3)$, $q =
(2, 2)$, $S = (\{ 2, 4, 6, 8 \}, \{ 2 \})$ and $T = (\{ 4, 6 \},
\varnothing)$, then $A$ is the path algebra of the quiver
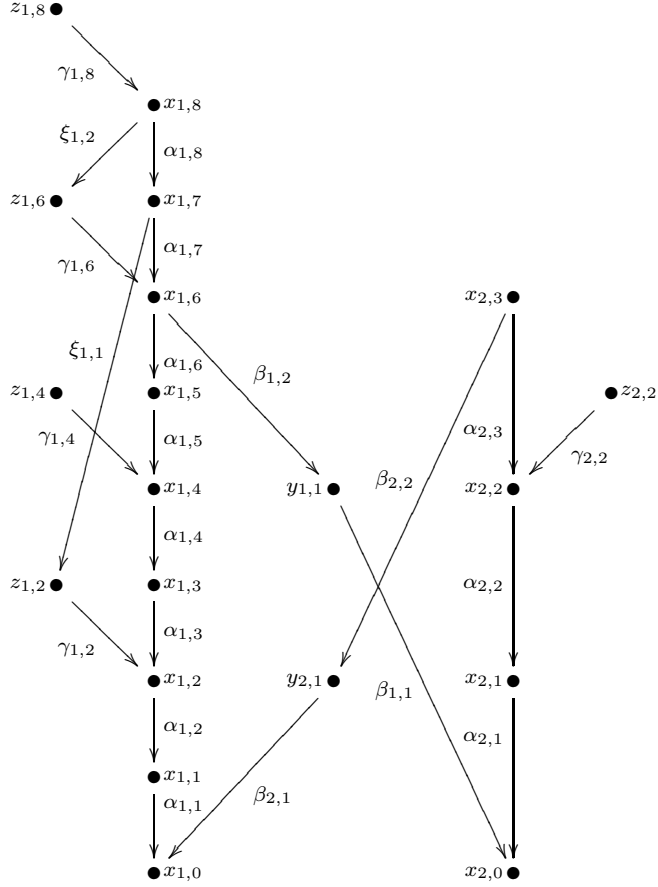
\[
\xymatrix{%
\bullet \save*+!R{\scriptstyle z_{1, 8}} \restore
\ar[rd]_{\gamma_{1, 8}} \\ %
& \bullet \save*+!L{\scriptstyle x_{1, 8}} \restore
\ar[d]^{\alpha_{1, 8}} \ar[ld]_{\xi_{1, 2}} \\ %
\bullet \save*+!R{\scriptstyle z_{1, 6}} \restore
\ar[rd]_{\gamma_{1, 6}} & \bullet \save*+!L{\scriptstyle x_{1, 7}}
\restore \ar[d]^{\alpha_{1, 7}} \ar[ldddd]_(.4){\xi_{1, 1}} \\ %
& \bullet \save*+!L{\scriptstyle x_{1, 6}} \restore
\ar[d]^(.7){\alpha_{1, 6}} \ar[rrdd]^{\beta_{1, 2}} & & & &
\bullet \save*+!R{\scriptstyle x_{2, 3}} \restore
\ar[dd]_(.7){\alpha_{2, 3}} \ar[lldddd]_{\beta_{2, 2}} \\ %
\bullet \save*+!R{\scriptstyle z_{1, 4}} \restore
\ar[rd]_(.3){\gamma_{1, 4}} & \bullet \save*+!L{\scriptstyle x_{1,
5}} \restore \ar[d]^{\alpha_{1, 5}} & & & & & \bullet
\save*+!L{\scriptstyle z_{2, 2}} \restore \ar[ld]^{\gamma_{2, 2}} \\ %
& \bullet \save*+!L{\scriptstyle x_{1, 4}} \restore
\ar[d]^{\alpha_{1, 4}} & & \bullet \save*+!R{\scriptstyle y_{1,
1}} \restore \ar[rrdddd]_{\beta_{1, 1}} & & \bullet
\save*+!R{\scriptstyle x_{2, 2}} \restore \ar[dd]_{\alpha_{2, 2}} \\ %
\bullet \save*+!R{\scriptstyle z_{1, 2}} \restore
\ar[rd]_{\gamma_{1, 2}} & \bullet \save*+!L{\scriptstyle x_{1, 3}}
\restore \ar[d]^{\alpha_{1, 3}} \\ %
& \bullet \save*+!L{\scriptstyle x_{1, 2}} \restore
\ar[d]^{\alpha_{1, 2}} & & \bullet \save*+!R{\scriptstyle y_{2,
1}} \restore \ar[lldd]^{\beta_{2, 1}} & & \bullet
\save*+!R{\scriptstyle x_{2, 1}} \restore \ar[dd]_(.3){\alpha_{2,
1}} \\ %
& \bullet \save*+!L{\scriptstyle x_{1, 1}} \restore
\ar[d]^(.3){\alpha_{1, 1}} \\ %
& \bullet \save*+!L{\scriptstyle x_{1, 0}} \restore & & & &
\bullet \save*+!R{\scriptstyle x_{2, 0}} \restore}
\]
bounded by relations
\begin{gather*}
\alpha_{1, 1} \alpha_{1, 2} \gamma_{1, 2}, \; \alpha_{1, 3}
\alpha_{1, 4} \gamma_{1, 4}, \; \alpha_{1, 5} \alpha_{1, 6}
\gamma_{1, 6}, \; \alpha_{1, 7} \alpha_{1, 8} \gamma_{1, 8}, \;
\alpha_{2, 1} \alpha_{2, 2} \gamma_{2, 2}, \; \beta_{1, 2}
\alpha_{1, 7},
\\ %
\xi_{1, 1} \alpha_{1, 8}, \; \alpha_{1, 2} \alpha_{1, 3}
\alpha_{1, 4} \alpha_{1, 5} \alpha_{1, 6} \alpha_{1, 7} -
\alpha_{1, 2} \gamma_{1, 2} \xi_{1, 1}, \; \alpha_{1, 6}
\alpha_{1, 7} \alpha_{1, 8} - \alpha_{1, 6} \gamma_{1, 6} \xi_{1,
2}.
\end{gather*}

In~\cite{BoSk1} we introduced a class of algebras called
admissible algebras with formal two-ray modules (see
also~\ref{subsectadmalg}). The main result of the paper is the
following.

\begin{theo} \label{theomain}
An algebra $A$ is an admissible algebra with formal two-ray
modules if and only if $A$ is the algebra of a defining system.
\end{theo}

Following~\cite{BoSk1} a connected translation quiver is said to
be of the first type, if its stable part is $\bbZ \bbA_\infty$,
while its left and right stable parts are $(-\bbN \bbD_\infty)$
and $\bbN \bbD_\infty$, respectively. Similarly, we say that a
connected translation quiver is of the second type, if its stable
part is the disjoint union of two quivers of the form $\bbZ
\bbA_\infty$, its left stable part is $(-\bbN) \bbA_\infty^\infty$
and its right stable part is the disjoint union of two quivers of
the form $\bbN \bbD_\infty$. As a consequence of the above theorem
and the Main Theorem of~\cite{BoSk1} we obtain the following.

\begin{theo} \label{theonext}
Let $A$ be the algebra of a defining system $(p, q, S, T)$ such
that $\sum_{i = 1}^{|p|} |S_i| > 0$. Let $L$ be the number of
indices $i \in [1, |p|]$ such that $|S_i| > 0$ and $\max S_i \in
T_i$. Then the Auslander--Reiten quiver of $A$ consists of the
following components:
\renewcommand{\labelenumi}{\textup{(\theenumi)}}
\begin{enumerate}

\item
a preprojective component of type $\tilde{\bbA}_{\sum_{i =
1}^{|p|} p_i, \sum_{i = 1}^{|p|} q_i}$,

\item
$\sum_{i = 1}^{|p|} |T_i| + 1$ families of coray tubes indexed by
$K$,

\item
$\sum_{i = 1}^{|p|} (|S_i| - |T_i|)$ components of the first type,

\item
$\sum_{i = 1}^{|p|} |T_i|$ components of the second type,

\item
a preinjective component of type $\tilde{\bbA}_{2, p_i + |T_i| -
\max S_i}$, for each $i \in [1, |p|]$ such that $|S_i| > 0$ and
$\max S_i \in T_i$,

\item
countably many components of the form $\bbZ \bbD_\infty$, if
$\sum_{i = 1}^{|p|} |T_i| > 0$,

\item
countably many components of the form $\bbZ \bbA_\infty^\infty$,
if $\sum_{i = 1}^{|p|} |T_i| > L$.

\end{enumerate}
\end{theo}

\section{Combinatorial structures} \label{sectcomb}

In this section we associate a combinatorial structure for two-ray
modules with a defining system.

\subsection{Definition} \label{subsectdef}

In order to define a notion of a combinatorial structure for
two-ray modules, we need a formalism of partial functions. Recall
that, if $I$ is a set, then every function of the form $\varphi :
D \to I$, where $D \subset I$, is called a partial function. We
write this fact $\varphi : I \dashrightarrow I$. The set $D$ is
called the domain of $\varphi$ and denoted $\Dom \varphi$. The
image of $\varphi$ will be denoted $\Im \varphi$. If $\varphi_1,
\varphi_2 : I \dashrightarrow I$, then the composition $\varphi_2
\varphi_1$ is the partial function $\varphi : I \dashrightarrow I$
with the domain $\varphi_1^{-1} (\Dom \varphi_2)$ defined by the
usual formula. If $\varphi : I \dashrightarrow I$ is injective,
then there exists a unique $\psi : I \dashrightarrow I$ such that
$\Dom \psi = \Im \varphi$ and $\psi \varphi = \Id_{\Dom \varphi}$.
In the above situation, $\Im \psi = \Dom \varphi$ and $\varphi
\psi = \Id_{\Im \varphi}$, and we call $\varphi$ an invertible
partial map and denote $\psi$ by $\varphi^-$. If $\varphi : I
\dashrightarrow I$ is arbitrary and $n \in \bbN$, then $\varphi^n$
denotes the $n$-fold composition of $\varphi$ with itself.
Moreover, by $\varphi^0$ we mean the identity map $\Id_I : I
\rightarrow I$. Finally, by $\varnothing : I \dashrightarrow I$ we
denote the empty map ($\Dom \varnothing = \varnothing$).

A $5$-tuple
\[
\langle I, \phi, \rho, \psi, (l_x)_{x \in (\Dom \phi \cup \Dom
\rho) \setminus \Dom \psi} \rangle
\]
consisting of a finite set $I$, invertible partial maps $\phi,
\rho, \psi : I \dashrightarrow I$, and a sequence $(l_x)_{x \in
(\Dom \phi \cup \Dom \rho) \setminus \Dom \psi}$ of nonpositive
integers, is called a combinatorial structure for two-ray modules,
if the conditions (C1)--(C14) listed below are satisfied. The
first ten conditions describe relations between the domains and
the images of $\phi$, $\rho$ and $\psi$:{
\renewcommand{\theenumi}{C\arabic{enumi}}
\begin{enumerate}

\item
$I = \Im \phi \cup \Im \rho \cup \Im \psi$,

\item
$\Im \phi \cap \Im \rho = \varnothing$,

\item
$\Im \phi \cap \Im \psi = \varnothing$,

\item
$\Im \rho \cap \Im \psi = \varnothing$,

\item
$\Dom \phi \cap \Dom \rho = \varnothing$,

\item
$\Dom \rho \cap \Im \rho = \varnothing$,

\item
$\Im \phi \subseteq \Dom \phi \cup \Dom \rho$,

\item
$\Dom \psi \subseteq \Dom \phi$,

\item
$\Im \rho \cap \Dom \psi = \varnothing$,

\item
$\Im \psi \subseteq \Dom \phi \cup \Dom \rho$.

\end{enumerate}}

The next condition is the following:{
\renewcommand{\theenumi}{C\arabic{enumi}}
\begin{enumerate}
\setcounter{enumi}{10}

\item
$\psi^{|I|} = \varnothing$.

\end{enumerate}
The above condition allows us to introduce two new partial
functions $\sigma, \eta : I \dashrightarrow I$. By definition,
$\Dom \sigma = \Dom \phi$ and $\sigma x = \psi^{v_{\phi x}} \phi
x$, where for $x \in I$, $v_x$ is the maximal nonnegative integer
$v$ such that $\psi^v x$ is defined (i.e., $x \in \Dom \psi^v$).
Similarly, we put $\Dom \eta = \Dom \phi \cup \Dom \rho$, $\eta x
= \sigma x$ for $x \in \Dom \phi$ and $\eta x = \rho x$ for $x \in
\Dom \rho$.}

The last three conditions put constrains on the numbers $l_x$:{
\renewcommand{\theenumi}{C\arabic{enumi}}
\begin{enumerate}
\setcounter{enumi}{11}

\item
if $x \in \Dom \psi$, then $l_{\sigma x} < 0$,

\item
if $x \in \Dom \phi \cap \Im \rho$, then $l_x = 0$ and $l_{\sigma
x} = 0$,

\item
for each $x \in I$, there exists $u \geq 0$ such that either
$\eta^u x \in \Im \rho \setminus \Dom \phi$ or $\sum_{k = 1}^u
l_{\eta^k x} < 0$ (the latter condition means in particular that
for all $k \in [1, u]$, $l_{\eta^k x}$ is defined, i.e., $\eta^k x
\in (\Dom \phi \cup \Dom \rho) \setminus \Dom \psi$). Note that by
the empty sum we always mean $0$.

\end{enumerate}
If there exists $u > 0$ such that $\sum_{k = 1}^u l_{\eta^k x} <
0$, then we set $u_x$ to be the maximal nonnegative integer $u$
such that $\sum_{k = 1}^u l_{\eta^k x} = 0$. Otherwise, by $u_x$
we denote the minimal $u \geq 0$ such that $\eta^u x \in \Im \rho
\setminus \Dom \phi$.}

\subsection{Notation} \label{subsectnot}

Before we associate a combinatorial structure with a defining
system, we introduce an additional notation. The reader is
encouraged to check the definitions below with the example
presented in Section~\ref{sectmain}.

Let $(p, q, S, T)$ be a defining system and let $Q$ be the quiver
of $(p, q, S, T)$ as defined in Section~\ref{sectmain}. We define
the following sets:
\begin{align*}
\fx & = \{ x_{i, j} \mid i \in [1, |p|], \, j \in [0, p_i + |T_i|] \},
\\ %
\fx_0 & = \{ x_{i, 0} \mid i \in [1, \ldots, |p|] \},
\\ %
\fx_1 & = \{ x_{i, j} \mid i \in [1, |p|], \, j \in [1,
p_i - 1] \},
\\ %
\fx_2 & = \{ x_{i, p_i} \mid i \in [1, |p|] \},
\\ %
\fx_3 & = \{ x_{i, j} \mid i \in [1, |p|], \, j \in [p_i + 1,
p_i + |T_i| - 1] \},
\\ %
\fx_4 & = \{ x_{i, p_i + |T_i|} \mid i \in [1, |p|] \},
\\ %
\fz & = \{ z_{i, j} \mid i \in [1, |S|], \, j \in S_i \}.
\end{align*}
We also define invertible partial functions $\fP, \fR, \fS, \fT :
\fx \cup \fz \dashrightarrow \fx \cup \fz$, by
\begin{align*}
\Dom \fP & = \fx \setminus \fx_0,
\\ %
\fP x_{i, j} & = x_{i, j - 1}, \, i \in [1, |p|], \, j \in
[1, p_i + |T_i|],
\\ %
\Dom \fR & = \fx_2 \cup \fx_3 \cup \fx_4,
\\ %
\fR x_{i, j} & =
\begin{cases}
x_{i + 1, 0} & i \in [1, |p|], \, j = |p_i|,
\\ %
x_{i, T_{i, j}} & i \in [1, |p|], \, j \in [|p_i| + 1, |p_i| +
|T_i|],
\end{cases}
\\ %
\Dom \fS & = \fz,
\\ %
\fS z_{i, j} & = x_{i, j}, \, i \in [1, |p|], \, j \in
S_i,
\\ %
\Dom \fT & = (\fx_3 \cup \fx_4) \setminus \fx_2,
\\ %
\fT x_{i, p_i + j} & = z_{i, T_{i, j}}, \, i \in [1, |p|], \, j
\in [1, |T_i|],
\end{align*}
where as usual $x_{n + 1, 0} = x_{1, 0}$. Note that $\fS \fT = \fR
|_{(\fx_3 \cup \fx_4) \setminus \fx_2}$.

We can describe $\fP$, $\fR$, $\fS$ and $\fT$ in a more intuitive
way. We divide the arrows of $Q$ into two groups: the first group
consists of the arrows $\alpha_{i, j}$, $i \in [1, |p|]$, $j \in
[1, p_i + |T_i|]$, while the second group consists of the
remaining arrows. A path in $Q$ is said to be of the first kind if
it is a composition of arrows from the first group, and of the
second kind if it is a composition of arrows from the second
group. Note that, for two paths of the first kind starting at the
same vertex, one of them has to be a subpath of the other. Thus,
for a given vertex $x$, we have a natural order in the set of all
paths of the first kind starting at $x$. The same remark applies
to the paths of the second kind starting at a given vertex, and to
the paths of a fixed kind terminating at a given vertex.

Now we give the foretold interpretations of the functions $\fP$,
$\fR$, $\fS$ and $\fT$. The function $\fP$ associates with $x \in
\fx \setminus \fx_0$ the terminating vertex of the minimal
nontrivial path of the first kind starting at $x$. Similarly, the
function $\fS$ (respectively, $\fT$) associates with $x \in \fz$
(respectively, $x \in (\fx_3 \cup \fx_4) \setminus \fx_2$) the
terminating vertex of the minimal nontrivial path of the second
kind starting at $x$. Finally, $\fR$ associates with $x \in \fx_2
\cup \fx_3 \cup \fx_4$ the terminating vertex of the minimal
nontrivial path of the second kind starting at $x$ whose
terminating vertex belongs to $\fx$. We leave to the reader to
figure out similar interpretations of the inverse functions
$\fP^-$, $\fR^-$, $\fS^-$ and $\fT^-$ (they involve paths of the
first and the second kind terminating at a given vertex).

For $x \in \Im \fS$, $x = x_{i, j}$, we write $\gamma_x$ for
$\gamma_{i, j}$, and for $x \in (\fx_3 \cup \fx_4) \setminus
\fx_2$, $x = x_{i, p_i + j}$, we denote $\xi_{i, j}$ by $\xi_x$.
Moreover, for $x \in \fx$, $x = x_{i, j}$, let $h_x = p_i + |\{ k
\in [1, |T_i|] \mid T_{i, k} \leq j \}|$ and $\omega_x =
\alpha_{i, j + 1} \cdots \alpha_{i, h_x}$. In particular, if $h_x
= j$, then $\omega_x$ is the trivial path at $x$, which we also
denote by $x$. Note that $h_x = j$ if and only if $x \in \fx_4$,
since $T_{i, j} \leq p_i + 2 j - |T_i| - 1$.

For each $x \in \fx$, we denote by $\mu_x$ the maximal path of the
second kind starting at $x$. Note that $\mu_x$ can be defined by
the following inductive rule:
\[
\mu_x =
\begin{cases}
x & x \in \fx_0 \cup \fx_1,
\\ %
\beta_{i, 1} \cdots \beta_{i, q_i} & x \in \fx_2, \, x = x_{i,
p_i},
\\ %
\mu_{\fR x} \gamma_{\fR x} \xi_x & x \in (\fx_3 \cup \fx_4)
\setminus \fx_2.
\end{cases}
\]
Moreover, $t \mu_x = \fR^j x$, where $j$ is the maximal
nonnegative integer $i$ such that $\fR^i x$ is defined. Recall
that, for a path $\tau$ in $Q$, $s \tau$ and $t \tau$ denote the
starting and the terminating vertex of $\tau$, respectively.

Similarly, for each $x \in \fx$, by $\nu_x$ we denote the maximal
path of the second kind terminating at $x$. We have
\[
\nu_x =
\begin{cases}
x & x \in \fx \setminus (\fx_0 \cup \Im \fS),
\\ %
\mu_{\fR^- x} \nu_{\fR^- x} & x \in \fx_0,
\\ %
\gamma_x & x \in \Im \fS \setminus \Im \fS \fT,
\\ %
\gamma_x \xi_{\fR^- x} \nu_{\fR^- x} & x \in \Im \fS \fT,
\end{cases}
\]
and
\[
s \nu_x =
\begin{cases}
(\fR^-)^j x & (\fR^-)^j x \not \in \Im \fS,
\\ %
\fS^- (\fR^-)^j x & (\fR^-)^j x \in \Im \fS,
\end{cases}
\]
where $j$ is the maximal nonnegative integer $i$ such that
$(\fR^-)^i x$ is defined.

\subsection{The combinatorial structure of a defining
system} \label{subsectcom}

Assume that $(p, q, S, T)$ is a defining system and let $Q$ be its
quiver. We will use the notation introduced in the previous
subsection.

Let $I = \fx \setminus \fx_4 \cup \fz$. We define partial
invertible functions $\phi, \rho, \psi : I \dashrightarrow I$ by
\begin{align*}
\Dom \phi & = I \setminus (\Im \fP \fS \cup \Im \fT),
\\ %
\phi x & =
\begin{cases}
t \mu_{\fP^- x} & x \in \fx \setminus (\fx_4 \cup \Im \fP \fS),
\\ %
t \mu_{\fS x} & x \in \fz \setminus \Im \fT,
\end{cases}
\\ %
\rho^- & = \fP \fS,
\\ %
\psi^- & = \fR |_{\Dom \fR \setminus \fx_4}.
\end{align*}
For each $x \in I$, we define a nonpositive integer $l_x$ by the
formula
\[
l_x =
\begin{cases}
- q_{i - 1} & x = x_{i, 1} \text { for } i \in [1, |p|] \text{
such that } p_i > 1,
\\ %
- q_{i - 2} & x = x_{i, 0} \text { for } i \in [1, |p|] \text{
such that } p_{i - 1} = 1,
\\ %
-2 & x \in \Dom \fP \text{ and } \fP x \in \Im \fS \fT,
\\ %
-2 & x \in \Im \fP \cap \Im \fR \text{ and }\fP^- x = \fR^- x,
\\ %
0 & \text{otherwise},
\end{cases}
\]
where $q_0 = q_n$, $q_{-1} = q_{n - 1}$, $p_0 = p_n$ and $T_0 =
T_n$. Note that $p_i = 1$ implies $T_i = \varnothing$.

Our aim in this subsection is to show that the above defined
structure $\langle I, \phi, \rho, \psi, (l_x)_{(\Dom \phi \cup
\Dom \rho) \setminus \Dom \psi} \rangle$ is a combinatorial
structure for two-ray modules, called the combinatorial structure
of $(p, q, S, T)$, and investigate its properties.

First observe that
\begin{align*}
\Dom \phi & = \fx \setminus (\fx_4 \cup \Im \fP \fS) \cup \fz
\setminus \Im \fT,
\\ %
\Im \phi & = \fx_0 \cup \fx_1,
\\ %
\Dom \rho & = \Im \fP \fS,
\\ %
\Im \rho & = \fz,
\\ %
\Dom \psi & = \Im \fR \setminus \fR (\fx_4),
\\ %
\Im \psi & = (\fx_2 \cup \fx_3) \setminus \fx_4,
\end{align*}
which immediately implies that the conditions~(C1)--(C10)
from~\ref{subsectdef} are satisfied. The condition~(C11) also
follows easily. Let $v_x$, $x \in I$, $\sigma$ and $\eta$ have the
same meaning as in~\ref{subsectdef}. One easily checks that
\[
\psi^{v_x} x =
\begin{cases}
s \nu_x & x \in \fx \setminus \fx_4, \, s \nu_x \in \fx \setminus
\fx_4,
\\ %
\fR s \nu_x & x \in \fx \setminus \fx_4, \, s \nu_x \in \fx_4,
\\ %
\fS s \nu_x & x \in \fx \setminus \fx_4, \, s \nu_x \in \fz, \,
\fS s \nu_x \not \in \fx_4,
\\ %
\fR \fS s \nu_x & x \in \fx \setminus \fx_4, \, s \nu_x \in \fz,
\, \fS s \nu_x \in \fx_4,
\\ %
x & x \in \fz.
\end{cases}
\]
Since
\begin{equation} \label{eqsnu}
s \nu_{\phi x} =
\begin{cases}
\fP^- x & x \in \fx \setminus (\fx_4 \cup \Im \fP \fS),
\\ %
x & x \in \fz \setminus \Im \fT,
\end{cases}
\end{equation}
we obtain
\begin{align}
\sigma x & =
\begin{cases}
\fP^- x & x \in \fx \setminus (\fx_4 \cup \Im \fP \fS), \, \fP^- x
\not \in \fx_4,
\\ %
\fR \fP^- x & x \in \fx \setminus (\fx_4 \cup \Im \fP \fS), \,
\fP^- x \in \fx_4,
\\ %
\fS x & x \in \fz \setminus \Im \fT, \fS x \not \in \fx_4,
\\ %
\fR \fS x & x \in \fz \setminus \Im \fT, \fS x \in \fx_4,
\end{cases} \label{eqsigma}
\\ %
\intertext{and} %
\notag \eta x & =
\begin{cases}
\fP^- x & x \in \fx \setminus (\fx_4 \cup \Im \fP \fS), \, \fP^- x
\not \in \fx_4,
\\ %
\fR \fP^- x & x \in \fx \setminus (\fx_4 \cup \Im \fP \fS), \,
\fP^- x \in \fx_4,
\\ %
\fS x & x \in \fz \setminus \Im \fT, \fS x \not \in \fx_4,
\\ %
\fR \fS x & x \in \fz \setminus \Im \fT, \fS x \in \fx_4,
\\ %
\fS^- \fP^- x & x \in \Im \fP \fS.
\end{cases}
\end{align}

It remains to verify the conditions~(C12)--(C14). A crucial
observation is that, if $x \in \Dom \phi$, then
\begin{equation} \label{eqlsigmax}
l_{\sigma x} < 0 \text{ if and only if } x \in \Im \fR,
\end{equation}
which follows by direct inspection. This immediately implies that
the condition~(C12) and the second part of the condition~(C13) are
satisfied. The first part of the condition~(C13) also follows
easily, since $\Dom \phi \cap \Im \rho = \fz \setminus \Im \fT$
and $l_x = 0$ for $x \in \fz$. Since $\Im \rho \setminus \Dom \rho
= \Im \fT$, in order to prove the condition~(C14), it is enough to
show that, for each $x \in I$, there exists $u \geq 0$ such that
$\eta^u x \in \Im \fT \cup \Im \fR$ (here again we
use~\eqref{eqlsigmax}).

Let $x \in I$. If $x \in \Im \fR \cup \Im \fT$, then the claim is
obvious. Assume now that $x \in \fx \setminus (\fx_4 \cup \Im
\fR)$. If $x \in \Im \fP \fS$, then $\eta x = \fS^- \fP^- x$. If
$\eta x \in \Im \fT$, then we are done, otherwise $\eta^2 x = \fR
\fP^- x \in \Im \fR$, if $\fP^- x \in \fx_4$, or $\eta^2 x = \fP^-
x$, if $\fP^- x \not \in \fx_4$. Similarly, if $x \not \in \Im \fP
\fS$, then $\eta x = \fR \fP^- x \in \Im \fR$, if $\fP^- x \in
\fx_4$, or $\eta x = \fP^- x$,  if $\fP^- x \not \in \fx_4$. Thus
the claim follows by easy induction for all $x \in \fx \setminus
\fx_4$. Finally, for $x \in \fz \setminus \Im \fT$, we have $\eta
x \in \fx \setminus \fx_4$, which finishes the proof.

The above considerations imply in particular that $u_x$, defined
as in~\ref{subsectdef}, is the minimal nonnegative integer $u$
such that $\eta^u x \in \Im \fR \cup \Im \fT$. It also follows
that, for $x \in \fx \setminus \fx_4$, $x = x_{i, j}$, $\eta^{u_x}
x \in \Im \fT$ if and only if $x \not \in \Im \fR$ and $h_x < p_i
+ |T_i|$. On the other hand, if $x \not \in \Im \fR$ and $h_x =
p_i + |T_i|$, then $\eta^{u_x} x = \fR x_{i, p_i + |T_i|}$.

\subsection{Admissible indices}

Let $\langle I, \phi, \rho, \psi, (l_x)_{x \in (\Dom \phi \cup
\Dom \rho) \setminus \Dom \psi} \rangle$ be a combinatorial
structure for two-ray modules. An index $y \in I$ is called
admissible if the following conditions are satisfied:{
\renewcommand{\theenumi}{A\arabic{enumi}}
\begin{enumerate}

\item
$y \in \Dom \phi$,

\item
$\sigma y \in \Dom \phi$,

\item
$l_{\sigma y} = 0$,

\item
if $y \in \Im \rho$, then $\eta^{u_y} y \in \Dom \phi$,

\item
$y \not \in \sigma (\Im \rho)$.

\end{enumerate}
The above definition is a modified version of the original
definition given in~\cite{BoSk1}. One may verify that the both
definitions are equivalent. The aim of this subsection is to
characterize admissible indices in the combinatorial structures of
defining systems.}

Let $(p, q, S, T)$ be a defining system. We use the notation
introduced in the previous two subsections.

\begin{lemm} \label{lemmadm}
Let $\langle I, \phi, \rho, \psi, (l_x)_{x \in (\Dom \phi \cup
\Dom \rho) \setminus \Dom \psi} \rangle$ be the combinatorial
structure for two-ray modules of $(p, q, S, T)$. Then $y \in I$ is
an admissible index if and only if either $y \in \fx \setminus
(\fx_0 \cup \fx_4 \cup \Im \fS \cup \Im \fP \fS \cup \Im \fP^2
\fS)$ or $y \in \fz \setminus \Im \fT$ and $h_{\fS y} = p_i +
|T_i|$ \textup{(}equivalently, $j > \max T_i$\textup{)}.
\end{lemm}

\begin{proof}
We first show that, if $y \in \fx$, $y = x_{i, j}$, and $y$ is
admissible, then $y \not \in \fx_0 \cup \fx_4 \cup \Im \fS \cup
\Im \fP \fS \cup \Im \fP^2 \fS$. Obviously, $y \not \in \fx_4$,
since $y \in I$. If $y \in \Im \fR$, then $l_{\sigma y} < 0$, and
hence $y$ is not admissible. If $y \in \Im \fS \setminus (\fx_4
\cup \Im \fR)$, then $y \in \sigma (\Im \rho)$, thus again $y$ is
not admissible. All together, this gives $y \not \in \fx_0 \cup
\fx_4 \cup \Im \fS$. Finally, if $y \in \Im \fP \fS$, then $y \not
\in \Dom \phi$, and if $y \in \Im \fP^2 \fS$, then $\sigma y \not
\in \Dom \phi$.

Now we check that, if $y \in \fx \setminus (\fx_0 \cup \fx_4 \cup
\Im \fS \cup \Im \fP \fS \cup \Im \fP^2 \fS)$, then $y$ is
admissible. First, observe that $y \in \Dom \phi$, since $y \in
\fx \setminus (\fx_4 \cup \Im \fP \fS)$. Moreover, if $\fP^- y
\not \in \fx_4$, then $\sigma y = \fP^- y \in \Dom \phi$, because
$y \not \in \Im \fP^2 \fS$. Otherwise, $\sigma y = \fR \fP^- y \in
\Dom \phi$, since $\Im \fR \subset \Dom \phi$. Next, using that $y
\not \in \Im \fR$ and~\eqref{eqlsigmax}, we get that $l_{\sigma y}
= 0$. The condition~(A4) is satisfied trivially, since $\Im \rho
\cap \fx = \varnothing$. The condition~(A5) also follows, because
$y \not \in \fx_0 \cup \Im \fS$.

Now we turn our attention to $y \in \fz$. If $y \in \Im \fT$, then
$y \not \in \Dom \phi$, and hence $y$ is not admissible. If $y
\not \in \Im \fT$ and $h_{\fS y} < p_i + |T_i|$, then in
particular $\fS y \not \in \fx_4$, so $\sigma y = \fS y$.
Consequently, $\eta^{u_y} y = \eta^{u_{\fS y}} \fS y \in \Im \fT$
according to the last remark in the previous subsection. Thus $y$
is not admissible, since $y \in \Im \rho$ and $\Im \fT \cap \Dom
\phi = \varnothing$. We leave to the reader to verify, that, if $y
\in \fz \setminus \Im \fT$ and $h_{\fS y} = p_i + |T_i|$, then $y$
is admissible.
\end{proof}

The last remark in the previous subsection implies that, if $y \in
\fz$, $y = z_{i, j}$, is an admissible index, then $\eta^{u_y} y =
\fR x_{i, p_i + |T_i|}$.

\subsection{Extensions of defining systems}

Let $(p, q, S, T)$ be a defining system and let $y$ be an
admissible index in the associated combinatorial structure. We
define families $S' = (S'_i)_{i = 1}^{|p|}$ and $T' = (T'_i)_{i =
1}^{|p|}$ by the following formulas. If $y  \in \fx$, $y = x_{i_0,
j_0}$, then $S_i' = S_i$ for $i \neq i_0$, $S_{i_0}' = S_{i_0}
\cup \{ j_0 + 1 \}$ and $T' = T$. If $y \in \fz$, $y = z_{i_0,
j_0}$, then $S' = S$, $T_i' = T_i$ for $i \neq i_0$ and $T_{i_0}'
= T_{i_0} \cup \{ j_0 \}$. It is easily seen that $(p, q, S', T')$
is a defining system, which we call the defining system obtained
by extension by $y$. Note that the quiver $Q$ of $(p, q, S, T)$ is
a subquiver of the quiver $Q'$ of $(p, q, S', T')$. Indeed, if $y
\in \fx$, $y = x_{i_0, j_0}$, then $Q'$ is obtained from $Q$ be
adding the vertex $x_{i_0, j_0 + 1}$ and the arrow $\gamma_{i_0,
j_0 + 1}$. If $y \in \fz$, $y = z_{i_0, j_0}$, then we add the
vertex $x_{i_0, p_{i_0} + |T_{i_0}| + 1}$ and the arrows
$\alpha_{i_0, p_{i_0} + |T_{i_0}| + 1}$ and $\xi_{i_0, |T_{i_0}| +
1}$.

We use for $(p, q, S, T)$ the notation introduced
in~\ref{subsectnot}. The analogous objects defined for $(p, q, S',
T')$ will be denoted by the same letter with $'$. We want to
describe in this subsection connections between objects defined
for $(p, q, S, T)$ and $(p, q, S', T')$.

First assume that $y \in \fx$, $y = x_{i_0, j_0}$. We have the
following easily verified formulas
\begin{gather*}
\fx_0' = \fx_0, \; \fx_1' = \fx_1, \; \fx_2' = \fx_2, \; \fx_3' =
\fx_3, \; \fx_4' = \fx_4, \; \fz' = \fz \cup \{ z_{i_0, j_0 + 1}
\},
\\ %
\fP' = \fP, \; \fR' = \fR, \; \fT' = \fT, \; \fS' x = \fS x, \, x
\in \fz, \; \fS' z_{i_0, j_0 + 1} = x_{i_0, j_0 + 1}.
\end{gather*}
Moreover, for $x \in \fx$,
\begin{gather} \label{eqpathsx}
\omega_x'  = \omega_x, \; \mu_x' = \mu_x, \; \nu_x' =
\begin{cases}
\nu_x \gamma_{i_0, j_0 + 1} & s \nu_x = x_{i_0, j_0 + 1},
\\ %
\nu_x & \text{otherwise}.
\end{cases}
\end{gather}
Similarly, if $y \in \fz$, $y = z_{i_0, j_0}$, then
\begin{gather*}
\fx_0' = \fx_0, \; \fx_1' = \fx_1, \; \fx_2' = \fx_2, \; \fx_3' =
\fx_3 \cup \{ x_{i_0, p_{i_0} + |T_{i_0}|} \},
\\ %
\fx_4' = \fx_4 \setminus \{ x_{i_0, p_{i_0} + |T_{i_0}|} \} \cup
\{ x_{i_0, p_{i_0} + |T_{i_0}| + 1} \}, \; \fz' = \fz,
\\ %
\fP' x = \fP x, \, x \in \fx \setminus \fx_0, \; \fP' x_{i_0,
p_{i_0} + |T_{i_0}| + 1} = x_{i_0, p_{i_0} + |T_{i_0}|},
\\ %
\fR' x = \fR x, \, x \in \fx_2 \cup \fx_3 \cup \fx_4, \; \fR'
x_{i_0, p_{i_0} + |T_{i_0}| + 1} = x_{i_0, j_0},
\\ %
\fT' x = \fT x, \, x \in (\fx_3 \cup \fx_4) \setminus \fx_2, \;
\fT' x_{i_0, p_{i_0} + |T_{i_0}| + 1}  = y, \; \fS'= \fS,
\end{gather*}
and, for $x \in \fx$, $x = x_{i, j}$,
\begin{gather} \label{eqpathsza}
\omega_x'  =
\begin{cases}
\omega_x \alpha_{i_0, p_{i_0} + |T_{i_0}| + 1} & i = i_0 \text{
and } j \geq j_0,
\\ %
\omega_x & \text{otherwise},
\end{cases}
\\ %
\mu_x' = \mu_x, \; \nu_x' =
\begin{cases}
\nu_x \xi_{i_0, |T_{i_0}| + 1} & s \nu_x = y,
\\ %
\nu_x & \text{otherwise}.
\end{cases} \label{eqpathszb}
\end{gather}
Finally,
\begin{gather*}
\omega_{x_{i_0, p_{i_0} + |T_{i_0}| + 1}} = x_{i_0, p_{i_0} +
|T_{i_0}| + 1},
\\ %
\mu_{x_{i_0, p_{i_0} + |T_{i_0}| + 1}} = \mu_{x_{i_0, j_0}}
\gamma_{i_0, j_0} \xi_{i_0, |T_{i_0}| + 1}, \; \nu_{x_{i_0,
p_{i_0} + |T_{i_0}| + 1}} = x_{i_0, p_{i_0} + |T_{i_0}| + 1}.
\end{gather*}

\subsection{Extensions of combinatorial structures}

Let
\[
\langle I, \phi, \rho, \psi, (l_x)_{x \in (\Dom \phi \cup \Dom
\rho) \setminus \Dom \psi} \rangle
\]
be a combinatorial structure for two-ray modules and let $y$ be an
admissible index. We recall from~\cite{BoSk1} the definition of
the combinatorial structure
\[
\langle I', \phi', \rho', \psi', (l_x')_{x \in (\Dom \phi' \cup
\Dom \rho') \setminus \Dom \psi'} \rangle
\]
obtained by extension by $y$. Choose an element $y'$ not in $I$.
We put $I' = I \cup \{ y' \}$, $\Dom \phi' = \Dom \phi \setminus
\{ y \} \cup \{ y' \}$, $\phi' x = \phi x$ for $x \in \Dom \phi
\setminus \{ y \}$ and $\phi' y' = \phi y$. In order to define the
remaining elements, we have to consider two cases.

Assume first that $y \not \in \Im \rho$. We put $\Dom \rho' = \Dom
\rho \cup \{ y \}$, $\rho' x = \rho x$ for $x \in \Dom \rho$,
$\rho' y = y'$ and $\psi'  = \psi$. Note that $(\Dom \phi' \cup
\Dom \rho') \setminus \Dom \psi' = (\Dom \phi \cup \Dom \rho)
\setminus \Dom \psi \cup \{ y' \}$. We define $l_x' = l_x$ for $x
\in (\Dom \phi \cup \Dom \rho) \setminus \Dom \psi$ and $l_{y'}' =
0$.

Assume now that $y \in \Im \rho$. Let $z = \eta^{u_{y}} y$.
Observe that $z \neq y$ and $z \in \Dom \phi \setminus \Dom \psi$.
We put $\rho' = \rho$, $\Dom \psi' = \Dom \psi \cup \{ z \}$,
$\psi' x = \psi x$ for $x \in \Dom \psi$ and $\psi' z = y'$. Note
that $(\Dom \phi' \cup \Dom \rho') \setminus \Dom {\psi'} = (\Dom
\phi \cup \Dom \rho) \setminus (\Dom \psi \cup \{ y, z \}) \cup \{
y' \}$. We define
\[
l_x' =
\begin{cases}
l_x & x \neq \sigma^2 y \text{ or } x = \sigma^2 y \text{ and } z
= \sigma y,
\\ %
-2 & x = y' \text{ and } z = \sigma y, \sigma^2 y \text{ or } x =
\sigma^2 y \text{ and } z \neq \sigma y,
\\ %
0 & x = y' \text{ and } z \neq \sigma y, \sigma^2 y.
\end{cases}
\]

We have the following correspondence between the above defined
extensions and the extensions of the defining systems.

\begin{lemm} \label{lemmext}
Let $(p, q, S, T)$ be a defining system, let $y$ be an admissible
index in the associated combinatorial structure, and let $(p, q,
S', T')$ be the defining system obtained by extension by $y$. Then
the combinatorial structure of $(p, q, S', T')$ is the same as the
combinatorial structure obtained from the combinatorial structure
of $(p, q, S, T)$ by extension by $y$.
\end{lemm}

\begin{proof}
The formulas for the functions follow directly by applying the
appropriate formulas listed in the previous subsection (one also
uses that, if $y \in \fz$, $y = z_{i_0, j_0}$, then $\eta^{u_y} y
= \fR x_{i_0, p_{i_0} + |T_{i_0}|})$. The formulas for the numbers
$l_x'$ also follow by case by case analysis, which is quite
tedious if $y \in \fz$, hence we omit it here.
\end{proof}

\section{Admissible algebras} \label{sectproof}

Throughout this section $(p, q, S, T)$ will be a fixed defining
system. We will use freely the notation introduced in the previous
section.

\subsection{Algebras with formal two-ray modules} \label{subsectadmalg}

We recall first the notion of a one-point extension of an algebra.
Let $A$ be an algebra and let $R$ be an $A$-module. By a one-point
extension of $A$ by $R$ we mean the matrix algebra of the form
\[
A [R] =
\begin{bmatrix}
A & R \\ 0 & K
\end{bmatrix}.
\]
Every $A$-module can be viewed also as an $A [R]$-module in the
obvious way. Moreover, if $X$ is an $A$-module, then by $\ol{X}$
we denote the $A [R]$-module defined on $X \oplus \Hom_A (R, X)$
by
\[
\begin{bmatrix}
a & r \\ 0 & \lambda
\end{bmatrix}
\begin{bmatrix}
x \\ f
\end{bmatrix}
=
\begin{bmatrix}
a x + f (r) \\ \lambda f
\end{bmatrix}.
\]
We refer to~\cite{Ri2}*{2.5} for more information on the one-point
extensions of algebras.

Recall from~\cite{BoSk1} that by an algebra with formal two-ray
modules we mean an algebra $A$ together with a combinatorial
structure for two-ray modules
\[
\langle I, \phi, \rho, \psi, (l_x)_{x \in (\Dom \phi \cup \Dom
\rho) \setminus \Dom \psi} \rangle
\]
and two collections $(X_i)_{i \in I}$, $(R_i)_{i \in \Dom \phi}$
of $A$-modules. If $y$ is an admissible index in the combinatorial
structure, then we may define a new algebra with formal two-ray
modules, called the algebra obtained by extension by $y$, in the
following way. We take $A' = A [R_y]$, the combinatorial structure
$\langle I', \phi', \rho', \psi', (l_x')_{x \in (\Dom \phi' \cup
\Dom \rho') \setminus \Dom \psi'} \rangle$ obtained by extension
by $y$, $X_x' = \ol{X}_x$ for $x \in I$, and $R_x' = \ol{R}_x$ for
$x \in \Dom \phi \setminus \{ y \}$. Finally, we put $X_{y'}' =
\tau_{A'} X_{\phi y}'$ and define $R_{y'}'$ to be the direct sum
of the middle terms of the Auslander--Reiten sequences starting at
the indecomposable direct summands of $X_{y'}'$.

We associate now with $(p, q, S, T)$ an algebra with formal
two-ray modules in the following way. Let $A$ be the algebra of
$(p, q, S, T)$ and let $\langle I, \phi, \rho, \psi, (l_x)_{x \in
(\Dom \phi \cup \Dom \rho) \setminus \Dom \psi} \rangle$ be the
combinatorial structure of $(p, q, S, T)$. For $x \in \fx
\setminus \fx_4$, let $X_x = M (\nu_x)$, and for $x \in \fz$, let
$X_x = M (\omega_{\fP \fS x})$. Here, for a path $\tau$ in $Q$, by
$M (\tau)$ we denote the corresponding string module (see for
example~\cite{BuRi}). Note that $\End_A (X_x) = k$ for all $x \in
I$.

In order to define the modules $R_x$, we need more information
about the modules $X_x$. We refer to \cite{Ri2}*{2.4} for details
about the method applied below in order to calculate the
Auslander--Reiten translation.

\begin{lemm} \label{lemmtau}
If $x \in \Dom \phi$, then $\tau_A X_{\phi x} = X_x$.
\end{lemm}

\begin{proof}
Assume first that $x \in \fx \setminus (\fx_4 \cup \Im \fP \fS)$.
Then $\nu_{\phi x} = \mu_{\fP^- x}$ and
\[
P (x) \xrightarrow{f} P (\fP^- x) \rightarrow X_{\phi x}
\rightarrow 0
\]
is a minimal projective presentation of $X_{\phi x}$. Using the
formula $\tau_A X_{\phi x} \simeq \Du \Coker \Hom_A (f, A)$, where
$\Du = \Hom_k (-, k)$ is the standard duality, we obtain that
$\tau_A X_{\phi x} \simeq M (\nu_x) = X_x$. Similarly, if $x \in
\fz \setminus \Im \fT$, then $\nu_{\phi x} = \mu_{\fS x}
\gamma_{\fS x}$ and
\[
P_{\fP \fS x} \rightarrow P_x \rightarrow X_{\phi x} \rightarrow 0
\]
is a minimal projective presentation of $X_{\phi x}$, hence we
obtain $\tau_A X_{\phi x} \simeq M (\omega_{\fP \fS x}) = X_x$.
\end{proof}

A consequence of the above lemma and the Auslander--Reiten formula
is that $\dim_K \Ext_A^1 (X_{\phi x}, X_x) = 1$, thus there is a
unique extension $R_x$ of $X_{\phi x}$ by $X_x$, which is not
isomorphic to $X_x \oplus X_{\phi x}$. The corresponding exact
sequence is the Auslander--Reiten sequence.

Assume, for a moment, that the defining system $(p, q, S, T)$ is
fundamental, i.e.\ $S_i = \varnothing = T_i$ for all $i \in [1,
|p|]$. In this case, the algebra with formal two-ray modules
associated with $(p, q, S, T)$ is also called fundamental. Note
that in this situation $A$ is a hereditary algebra of type
$\tilde{\bbA}_{\sum_{i = 1}^{|p|} p_i, \sum_{i = 1}^{|q|} q_i}$,
the elements of $I$ are in bijection with the rays in a chosen
nonhomogeneous tube in the Auslander--Reiten quiver of $A$,
$\phi^-$ corresponds to the action of $\tau_A$, $\rho =
\varnothing = \psi$, $X_x$ are the corresponding simple regular
modules, $R_x$ are the corresponding modules of regular length
$2$, and $l_x = 1 - \dim_K \tau_A X_x$. An algebra with formal
two-ray modules is called admissible if it can be obtained from a
fundamental one by a sequence of extensions by admissible indices.

\subsection{Homomorphisms between modules}

We want to show that the extensions of the defining systems
correspond to the extensions of the associated algebras with
formal two-ray modules. In order to do this, we need a more
precise knowledge about the homomorphism spaces between the
corresponding modules. As above, $(p, q, S, T)$ is a fixed
defining system. We will use the notation introduced in the
previous subsection for the algebra with formal two-ray modules
associated with $(p, q, S, T)$.

We start with a remark about the homomorphism spaces between
string modules. Let $\tau_1$ and $\tau_2$ be paths in the quiver
$Q$ of $(p, q, S, T)$. Then it follows by easy calculations
(compare also~\cites{CB1, Sc}), that we have $\dim_K \Hom_A (M
(\tau_1), M (\tau_2)) = 1$ if and only if $\tau_1 = \tau_1'
\tau_0$ and $\tau_2 = \tau_0 \tau_2'$ for paths $\tau_0$,
$\tau_1'$ and $\tau_2'$ in $Q$, and $\Hom_A (M (\tau_1), M
(\tau_2)) = 0$ otherwise. A direct consequence of the above
formula is the following.

\begin{lemm} \label{lemmhomx} \rule{0pt}{0pt}
\renewcommand{\labelenumi}{\textup{(\theenumi)}}
\begin{enumerate}

\item
Let $y \in \fx \setminus \fx_4$ be such that $s \nu_y \not \in \Im
\fP \fS$ and let $x \in I$. Then $\dim_K \Hom_A (X_y, X_x) = 1$ if
$x \in \fx$ and $x = \fR^{-k} y$ for some $k \geq 0$, and $\Hom_A
(X_y, X_x) = 0$ otherwise.

\item
Let $y \in \fz \setminus \Im \fT$, $y = z_{i_0, j_0}$, be such
that $h_y = p_i + |T_i|$ and let $x \in I$. Then $\dim_K \Hom_A
(X_y, X_x) = 1$ if $x \in \fz$, $x = z_{i_0, j}$ for $j \in [j_0 +
1, p_i + |T_i|]$, and $\Hom_A (X_y, X_x) = 0$ otherwise. \qed

\end{enumerate}
\end{lemm}

The next step is the following.

\begin{lemm} \label{lemmhomrx}
Let $y \in \Dom \phi$ be an admissible index and $x \in I$. Then
$\dim_K \Hom_A (R_y, X_x) = 1$ if one of the following conditions
is satisfied:
\begin{itemize}

\item
$y \in \fx$, $y = x_{i_0, j_0}$, $x \in \fx$ and $s \nu_x =
x_{i_0, j_0 + 1}$,

\item
$y \in \fz$, $x \in \fx$ and $s \nu_x = y$,

\item
$y \in \fz$, $y = z_{i_0, j_0}$, and $x = z_{i_0, j}$ for $j \in
[j_0 + 1, p_{i_0} + |T_{i_0}|]$,

\end{itemize}
and $\Hom_A (R_x, X_y) = 0$ otherwise.
\end{lemm}

\begin{proof}
Since
\[
0 \rightarrow X_y \rightarrow R_y \rightarrow X_{\phi y}
\rightarrow 0
\]
is an Auslander--Reiten sequence, thus applying $\Hom_A (-, X_x)$
we get a short exact sequence
\[
0 \rightarrow \Hom_A (X_{\phi y}, X_x) \rightarrow \Hom_A (R_y,
X_x) \rightarrow \rad_A (X_y, X_x) \rightarrow 0.
\]

Assume first that $y \in \fx$. Then $s \nu_y = y \not \in \Im \fP
\fS$ and $s \nu_{\phi y} = \fP^- y \not \in \Im \fP \fS$, by
Lemma~\ref{lemmadm}. Moreover, $\rad_A (X_y, X_x) = 0$. Indeed,
from Lemma~\ref{lemmhomx}, the condition $\rad_A (X_y, X_x) \neq
0$ implies $x = \fR^{-k} y$ for some $k > 0$ (we exclude $k = 0$,
since $\rad_A (X_y, X_x) \neq 0$ implies that $y \neq x$), which
is impossible since $y \not \in \Im \fR$, by Lemma~\ref{lemmadm}.
Thus $\Hom_A (R_y, X_x) \simeq \Hom_A (X_{\phi y}, X_x)$ and,
using once more Lemma~\ref{lemmhomx}, we conclude that $\Hom_A
(R_y, X_x) \neq 0$ if and only if $x = \fR^{-k} \phi y$ for some
$k \geq 0$. One easily checks that this is equivalent to the first
condition in the lemma.

Now, let $y \in \fz$, $y = z_{i_0, j_0}$. Recall that, by
Lemma~\ref{lemmadm}, we have $y \not \in \Im \fT$ and $h_y =
p_{i_0} + |T_{i_0}|$. Moreover, $\phi y \in \fx$ and $s \nu_{\phi
y} = y \not \in \Im \fP \fS$. Hence, it follows immediately from
Lemma~\ref{lemmhomx} that either $\Hom_A (X_y, X_x) = 0$ or
$\Hom_A (X_{\phi y}, X_x) = 0$. If the first condition is
satisfied, then $\Hom_A (R_y, X_x) \simeq \Hom_A (X_{\phi y},
X_x)$, and consequently $\Hom_A (R_y, X_x) \neq 0$ if and only if
$x = \fR^{-k} \phi y$ for some $k \geq 0$. This leads to the
second condition in the lemma. In the latter case we have $\Hom_A
(R_y, X_x) \simeq \rad_A (X_y, X_x)$ and we get the third
condition (we exclude $j = j_0$, since then $y = x$ and $\rad_A
(X_y, X_x) = 0$).
\end{proof}

The final lemma is the following.

\begin{lemm} \label{lemmhomr}
If $y \in D_\phi$ is an admissible index and $x \in D_\phi$, $x
\neq y$, then $\Hom_A (R_y, X_{\phi x}) = 0$.
\end{lemm}

\begin{proof}
If $y \in \fx$, $y = x_{i_0, j_0}$, then $\Hom_A (R_y, X_{\phi x})
\neq 0$ implies $s \nu_{\phi x} = x_{i_0, j_0 + 1}$. If $x \in
\fx$, $x = x_{i, j}$, then, by~\eqref{eqsnu}, $s \nu_{\phi x} =
x_{i, j + 1}$, and hence $x = y$, a contradiction to the
hypothesis. If $x \in \fz$, $x = z_{i, j}$, then applying
again~\eqref{eqsnu}, we get $s \nu_{\phi x} = x_{i, j}$. Hence $y
= \fP \fS x$, which is again impossible, since $y \not \in \Im \fP
\fS$.

Assume now that $y \in \fz$, $y = z_{i_0, j_0}$. Since always
$\phi x \in \fx$, thus $\Hom_A (R_y, X_{\phi x}) \neq 0$ implies
$s \nu_{\phi x} = y$. By~\eqref{eqsnu} this is possible only if $x
= y$, a contradiction.
\end{proof}

\subsection{Proofs}

Now we give proofs of the main results. In order to prove
Theorem~\ref{theomain}, it is enough to show the following two
claims. First, if $y$ is an admissible index in the combinatorial
structure of $(p, q, S, T)$ and $(p, q, S', T')$ is the defining
system obtained by extension by $y$, then the algebra with formal
two-ray modules associated with $(p, q, S', T')$ is the extension
by $y$ of the algebra with formal two-ray modules associated with
$(p, q, S, T)$. Second, the defining system $(p, q, S, T)$ can be
obtained by a sequence of extensions by admissible indices from a
fundamental defining system. The latter claim is an easy
observation, hence we will concentrate on the former one. For the
algebra with formal two-ray modules associated with $(p, q, S, T)$
we will use the notation introduced above, while for the algebra
with formal two-ray modules associated with $(p, q, S', T')$, we
will use the analogous notation with $'$.

By direct calculations it follows that $A' = A [R_y]$. The
required relationship between the combinatorial structures is
given in Lemma~\ref{lemmext}. The formulas $X_x' = \ol{X}_x$, for
$x \in I$, follow from Lemma~\ref{lemmhomrx} and~\eqref{eqpathsx},
\eqref{eqpathsza}, \eqref{eqpathszb}. Next, the formulas for
$X_{y'}'$ and $R_{y'}'$ are consequences of Lemma~\ref{lemmtau}
and the definition of $R_{y'}'$. It remains to show that
\[
0 \to \ol{X}_x \to \ol{R}_x \to \ol{X}_{\phi x} \to 0
\]
is an Auslander--Reiten sequence of $A'$-modules, for all $x \in
\Dom \phi$, $x \neq y$. Since $\ol{X}_{\phi x} = X_{\phi x}$,
according to Lemma~\ref{lemmhomr}, this is a consequence
of~\cite{Ri2}*{2.5(6)}

For the proof of Theorem~\ref{theonext}, recall from
\cite{BoSk1}*{Main Theorem} (and its proof), that the
Auslander--Reiten quiver of $A$ consists of the following
components:{
\renewcommand{\labelenumi}{\textup{(\theenumi)}}
\begin{enumerate}

\item
a preprojective component of type $\tilde{\bbA}_{\sum_{i =
1}^{|p|} p_i, \sum_{i = 1}^{|p|} q_i}$,

\item
$N + 1$ families of coray tubes indexed by $K$,

\item
$M - N$ components of the first type,

\item
$N$ components of the second type,

\item
a preinjective component of type $\tilde{\bbA}_{2, |c|}$, for each
$\sigma$-cycle $c$ ,

\item
countably many components of the form $\bbZ \bbD_\infty$, if $N >
0$,

\item
countably many components of the form $\bbZ \bbA_\infty^\infty$,
if $N > L$,

\end{enumerate}
where $M = |\Dom \rho|$, $N = |\Dom \psi|$ and $L$ is the number
of the $\sigma$-cycles. Now $|\Dom \rho| = |\Im \rho| = |\fz| =
\sum_{i = 1}^{|p|} |S_i|$ and $|\Dom \psi| = |\Im \psi| = |(\fx_2
\cup \fx_3) \setminus \fx_4| = |(\fx_2 \cup \fx_3 \cup \fx_4)
\setminus \fx_4| = |p| + \sum_{i = 1}^{|p|} |T_i| - |p| = \sum_{i
= 1}^{|p|} |T_i|$. Finally, using~\eqref{eqsigma}, one checks that
the $\sigma$-cycles are of the form $\{ x_{i, j}, x_{i, j + 1},
\ldots, x_{i, p_i + |T_i| - 1} \}$, for $i \in [1, |p|]$ and $j
\in T_i$ such that $l \not \in S_i$ for $l \in [j + 1, p_i +
|T_i|]$.}

\end{document}